\documentclass[11pt]{article}
\usepackage[english]{babel}
\usepackage{amsfonts}
\usepackage[dvips]{graphicx}

\begin{document}

\title{\bf Fast computing of velocity field for flows in industrial burners and pumps}

\date{June 2007}

\author{\bf Gianluca Argentini \\
\normalsize gianluca.argentini@gmail.com \\
\normalsize gianluca.argentini@riellogroup.com \\
\textit{Research \& Development Department}\\
\textit{Riello Burners}, 37045 San Pietro di Legnago (Verona), Italy}

\maketitle

\begin{abstract}
In this work we present a technique of fast numerical computation for solutions of Navier-Stokes equations in the case of flows of industrial interest. At first the partial differential equations are translated into a set of nonlinear ordinary differential equations using the geometrical shape of the domain where the flow is developing, then these ODEs are numerically resolved using a set of computations distributed among the available processors. We present some results from simulations on a parallel hardware architecture using native multithreads software and simulating a shared-memory or a distributed-memory environment.\\

\noindent \textbf{keywords}: Navier-Stokes equations, streamlines, nonlinear differential equations, parallel computation.
\end{abstract}

\section{Introduction}
\label{intro}

\noindent

Design, development and engineering of industrial power burners have strong mathematical requests, such as numerical resolution of differential problems involving {\it Navier-Stokes} equations for velocity and pressure fields, geometrical design of combustion heads for a correct shape and optimal efficiency of {\it flame}, geometrical design of {\it ventilation fans} and computation of correct air inflows for optimal combustion, computation of flows in the internal vanes of {\it oil pumps} for estimating the correct mass rate and dimensions of gears.\\
Rapid prototyping for an accurate design of the correct geometries involves a numerical simulation of the gas or oil flows in these burner's components.\\
The necessity of an high graphic resolution for the visualization of the flows requires the knowledge of a large amount of streamlines. Therefore the numerical computation is very onerous about amount of memory and power of processors of the used hardware environment. A technique for reducing the total time of computation can be very useful. In this paper we present a method based on reduction of Navier-Stokes system to a set of not coupled ordinary differential equations. Each of these equations describes the fluid velocity field along each flow streamline. The advantage of this method, compared to an usual one such as finite differences scheme (see e.g. \cite{quarteroni}), is a faster computation of the numerical resolution and an easy parallelization in a multicore environment, using a distribution of blocks of streamlines among the available processors.

\begin{figure}[ht]\label{combustionHead}
	\begin{center}
	\includegraphics[width=8cm]{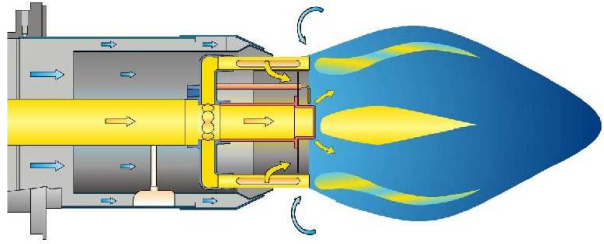}
	\caption{\small{\it Flows (gas at central zone) in a burner's combustion head. The swirling effects are very important for an optimized combustion, and they are very difficult to simulate by standard numerical methods.}} 
	\end{center}
\end{figure}

\section{The mathematical model}
\label{model}

In this section we briefly present the mathematical model developed for simbolic and numerical resolution of the differential problems about some flows in burner's components.\\
The main set of equations is the Navier-Stokes system describing the dynamics of the fluids (see e.g. \cite{anderson}):

\begin{equation}\label{NS}
	\rho \left( \frac{\partial \bf{v}}{\partial t} + \bf{v} \cdot \nabla \bf{v} \right) = - \nabla p + \mu \Delta \bf{v}
\end{equation}

\noindent where $\bf{v}$ is the velocity field, $p$ is the pressure, $\rho$ the density and $\mu$ the dynamic viscosity of the gas or oil. Sometimes there are other equations to be considered, such as continuity equation, but here we would describe a general simplified case. A geometric characterization of the flow is offered by the {\it streamlines}, that is the continuous lines on which the tangent at any point is parallel to the velocity vector.\\
\noindent At every instant of time, the velocity field in the domain of the flow is the set of vectors which give the velocities of the particles at every point ({\it eulerian description} of a flow). The computation of this vectorial field is done by numerical discretizations of (\ref{NS}), e.g. using finite temporal steps, which are usually expensive for CPUs and central memory of an hardware environment.\\
\noindent The fluid motion can be described by the set of all the trajectories followed by the particles and the speed along these geometric lines in a given interval of time ({\it lagrangian description} of a flow). The computation of the trajectories usally requires the knowledge of the velocity field, therefore numerical methods are expensive too (see \cite{chung} for a review on computational fluid dynamics techniques).\\
\noindent We present a method for computing the velocity field by reduction of Navier-Stokes equations to ordinary differential form.\\

Let $\Phi \in C^k([r_1,r_M] \times [s_0,s_1] \times [t_i,t_f], \mathbb{R}^3)$ be a smooth $(k \geq 1)$ function such that, if $t_0 \in [t_i,t_e]$ is a time value between initial $t_i$ and final instant $t_f$, and $r_{\phi} \in [r_1,r_M]$, the {\it curve} $\phi: s \longmapsto \Phi(r_{\phi},s,t_0)$ is a $C^k([s_0,s_1], \mathbb{R}^3)$ parameterized invertible representation (see \cite{lipschutz}) of a streamline at time $t_0$. The $r$ parameter identifies a single streamline, and in numerical computation it takes integer values in interval ${1,M}$, where $M$ is the total number of streamlines.
Then temporal evolution of the flow is described by the map $t \longmapsto \Phi(r,s,t)$, for $(r,s) \in [r_1,r_M]\times[s_0,s_1]$. At any $t_0 \in [t_i,t_f]$, the map $(r,s) \longmapsto \Phi(r,s,t_0)$ is a snapshot of all the streamlines at time $t_0$. We can consider the velocity field of this snapshot as a set of vectors depending only from $(r,s)$ and not from $t$, so that for the instant $t_0+dt$ the field is computed using the new snapshot $s \longmapsto \Phi(r,s,t_0+dt)$. \\
For fixed $t_0$ and $r_{\phi}$, let $(x_i)_{1 \leq i \leq 3}$ be a cartesian coordinates system in $\mathbb{R}^3$ and $s$ the parameter of the representation $s \longmapsto \phi(s)=\Phi(r_{\phi},s,t_0)$. In general, under usual hypothesis of regularity, the relation between a streamline and its tangent vector field is

\begin{equation}\label{tgRelation}
	\frac{v_i}{v_j} = \frac{\dot{\phi}_i}{\dot{\phi}_j}
\end{equation}

\noindent where $\dot{\phi}_i = \frac{d\phi_i}{ds}$. Now we define the parameterized velocity field {\it along} the streamline as 

\begin{equation}\label{paramField}
	\bf{u} = \bf{v} \circ \phi
\end{equation}

\noindent so that ${\bf u} = {\bf u}(r_{\phi},s,t_0) = {\bf u}(s)$. Using chain rule we have $\dot{u}_i = \sum_{m=1}^3 v_{im}\dot{\phi}_m$ and, from (\ref{tgRelation}), the following relation holds:

\begin{equation}
	\dot{u}_i = \left(v_1 v_{i1} + v_2 v_{i2} + v_3 v_{i3}\right)\frac{\dot{\phi}_i}{v_i}
\end{equation}

\noindent Note that the expression $v_1 v_{i1} + v_2 v_{i2} + v_3 v_{i3}$ is the $i$-th component of $\bf{v} \cdot \nabla \bf{v}$ in Navier-Stokes system, so it can be substituted by $v_i \dot{u}_i (\dot{\phi}_i)^{-1}$. \\
Now let $g = g({\bf x})$ be the inverse of $\phi$ and $q = p \circ \phi$; using standard computations we obtain

\begin{eqnarray}
	\frac{\partial v_i}{\partial x_j} & = & \dot{u}_i\frac{\partial g}{\partial x_j} \\
	\frac{\partial^2 v_i}{\partial x_j^2} & = & \ddot{u}_i\left(\frac{\partial g}{\partial x_j}\right)^2 + \dot{u}_i\frac{\partial^2 g}{\partial x_j^2} \\
	\frac{\partial p}{\partial x_j} & = & \dot{q}\frac{\partial g}{\partial x_j}
\end{eqnarray}

\noindent Along the streamline $s \longmapsto \phi(s)$ we have $\partial_t {\bf v} = \partial_t ({\bf u} \circ g) = {\bf 0}$, so Navier-Stokes equations can be translated into the following ordinary differential equations for the functions $(u_i)_{1 \leq i \leq 3}$:

\begin{equation}\label{odeNS}
	\rho u_i \dot{u}_i \dot{\phi}_i = - \dot{q}\frac{\partial g}{\partial x_i} + \mu \left( \nabla g \cdot \nabla g \hspace{0.1cm} \ddot{u}_i + \Delta g \hspace{0.1cm} \dot{u}_i \right)
\end{equation}

\noindent where $\rho$ and $\mu$ should be considered as function of the parameter $s$ and the derivatives of $g$ should be expressed throught those of $\phi$. These equations are again nonlinear, and a possible symbolic resolution required non standard analytical methods (as {\it Lie symmetries}, see \cite{stephani}), but now all the derivatives are ordinary and the equations are no more coupled.

\section{Example: internal gears pump}
\label{examplePump}

As example of technological interest, we consider an internal {\it gears pump} for oil, usually formed by a central internal rotor (pignon) which moves an external one (crown). 

\begin{figure}\label{lobePump}
\begin{center}
\includegraphics[width=6cm]{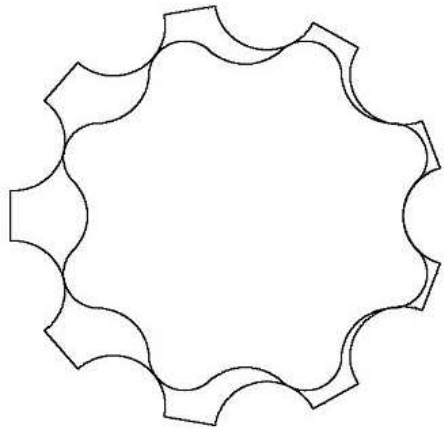}
\caption{\scriptsize{The particular case of {\it lobes} gear pump. The internal curve is the profile of the pignon, the external the profile of crown. The internal zones are full of oil, which is transported from suction inlet to delivery outlet. The analytical equations of the curves are of type $(x,y)=[r+sin(\phi+\omega t)]\left(cos(\pi+t),sin(\pi+t)\right)$.}}
\end{center}
\end{figure}

The oil is contained in the geometrical zones ({\it vanes}) between the two rotors. The knowledge of the flow in these vanes is important for a comprehension of some critical phenomena as cavitation (bubbles incoming) or loss of lubrication.\\
We would study the flow in the critical zone of oil compression, where the rotors are very near and the vane has a small volume. Here we consider, for simplicity, that streamlines have the same geometrical shape of the rotors. In the case, e.g., of a classical teeth gears pump, without considering a dimensional factor, a single profile (left side of the tooth) can be described by the following representation

\begin{eqnarray}\label{toothEqs}
	\phi_1 & = & cos(s) + s \hspace{0.1cm} sin(s) \\
	\nonumber \phi_2 & = & sin(s) - s \hspace{0.1cm} cos(s)
\end{eqnarray}

\noindent where $\pi + arctan(s_0) \leq s \leq \frac{3}{2}\pi$, for an opportune value $s_0$.

\begin{figure}\label{tooth}
\begin{center}
\includegraphics[width=6cm]{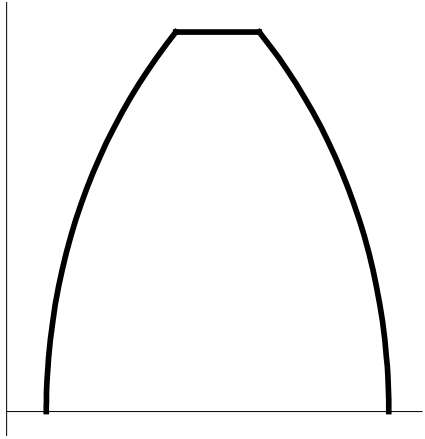}
\caption{\scriptsize{Tooth profile for internal rotor of a standard gear pump. The (\ref{toothEqs}) are the equations of the left side.}}
\end{center}
\end{figure}

The evolution in time can be handled, e.g., by a rotation of angle $\theta=\omega t$ for an interval of time small enough to have a costant shape for the volume between the two rotors. Then in this interval the streamlines have a new parametrization given by

\begin{equation}
 \Phi(s,t) = \phi_t(s) = \left[ \begin{array}{cc}
cos(\omega t) & \hspace{0.1cm} -sin(\omega t)  \\
sin(\omega t) & \hspace{0.1cm} cos(\omega t)  \end{array} \right] \phi(s)
\end{equation}

\noindent but the derivatives in the form (\ref{odeNS}) of Navier-Stokes equations are referred only to $s$ parameter again. With standard calculations it can be shown that, in the simple case of constant pressure, these equations have the nonlinear form

\begin{equation}
	\rho u_i \dot{u}_i = \mu \hspace{0.1cm} cos(s) \left[ \frac{1+s^2}{s}\ddot{u}_i + \frac{s^2-1}{s^2}\dot{u}_i \right]
\end{equation}

\section{The computational algorithm}
\label{complex}

\begin{figure}\label{suction}
\begin{center}
\includegraphics[width=7cm]{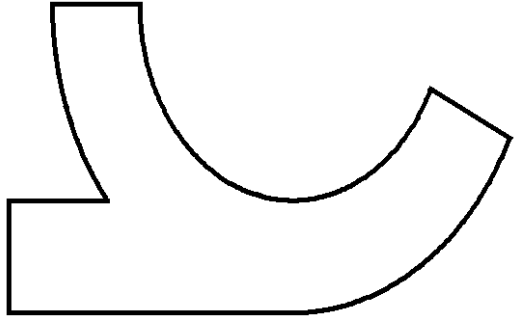}
\caption{\scriptsize{Suction zone for lobe pump. Oil flows from left horizontal channel and is distributed, by drop of pressure, among the lobes into the pump.}}
\end{center}
\end{figure}

For computing a solution of (\ref{odeNS}), also in case of more complex geometries than previous or in case of turbolent flow, we have tested an algorithm based on smooth interpolation of streamlines which are computed by a first rough numerical resolution of original Navier-Stokes equations. For example, in the oil suction channel (see fig. \ref{suction}) there is a geometric angle where the flow can become turbulent and problems of cavitation (formation of air or vapours bubbles due to local gradient of pressure) can happen.

\noindent We describe now the algorithm. For a chosen set of time values, a suitable number of streamlines are computed using a finite difference scheme on a grid covering the geometrical domain of interest (see \cite{anderson}). Then each streamline is divided into couples of geometrical consecutive points, i.e. $\{P_1,P_2\}$, $\{P_2,P_3\}$, ..., $\{P_{N-1},P_N\}$. Every couple is interpolated with a cubic polynomial $p=p(t)$, $t\in[0,1]$ ($t$ is not time), imposing the following four analytical conditions:

\begin{equation}
\left\{
\begin{array}{lll}\label{systemCubic}
	p(0)  &=& P_k\\
	p(1)  &=& P_{k+1}\\
	p'(0) &=& v_k\\
	p'(1) &=& v_{k+1}
\end{array}
\right.
\end{equation}

\noindent where $1 \leq {\it k} \leq N-1$ and $v$ is the numerical solution (velocity field) of the original Navier-Stokes system (see \cite{argentini} for details). In this way we obtain a set of class ${\it C}^1$ new streamlines. Note that cubics are the smallest degree one variable polinomials useful to have continuous differentiation along the interpoled lines. Each new streamline can be represented using, e.g., a unique parameter $s \in [0,1]$ with the re-parametrization $s=\frac{(m+t)}{N}$, with $m \in \mathbb{N}$, $0 \leq m \leq N-1$; if $f = f(t)$ is a cubic between two generic points, the curve $\phi = \phi(s)$ has derivative

\begin{equation}
	\dot\phi = N\frac{df}{dt}
\end{equation}

\noindent and two consecutive cubics have the same derivative in their common point. The first order derivative $\dot\phi$ is necessary for the ODE form (\ref{odeNS}) of the Navier-Stokes equations. It can be shown (see \cite{argentini}) that, if $M$ is the number of streamlines which we want to interpolate, all the $4 \times (N-1) \times M$ coefficients of the cubics can be computed by $N-1$ matrix-vector products ${\bf A}{\bf b}$ where ${\bf A}$ is the $4M \times 4M$ matrix\\

${\bf A} = \left( \begin{array}{ccccc}
{\bf T}  & {\bf 0} & . & . &  {\bf 0} \\
{\bf 0}  & {\bf T} & . & . &  {\bf 0} \\
. & . & . & . & . \\
. & . & . & . & . \\
{\bf 0}  & {\bf 0} & . & . &  {\bf T}
\end{array} \right)$,\\

\noindent ${\bf T}$ is the constant submatrix\\

${\bf T} = \left( \begin{array}{cccc}
2  & -2 & 1  & 1 \\
-3 & 3  & -2 & -1 \\
0  & 0  & 1  & 0 \\
1  & 0  & 0  & 0 
\end{array} \right)$\\

\noindent and ${\bf b}$ is the vector ${\bf b} =( p_{(1,k)}, p_{(1,k+1)}, . . ., v_{(M,k)}, v_{(M,k+1)})$, $1 \leq {\it k} \leq N-1$. Note that {\bf A} is a {\it sparse} matrix with {\it density number} (see \cite{golub}) at most $\frac{1}{M}$. Also, {\bf A} is a block diagonal matrix where each block is constant, so the execution of a product like ${\bf A}{\bf b}$ can be optimized using standard techniques (see \cite{golub}).

\section{About numerical resolution}
\label{numericalRes}

In this section we discuss the numerical resolution of the nonlinear differential equations (\ref{odeNS}).\\
For physical reasons, as typical in eulerian point of view, we might resolve an initial value problem for every streamline and for every instant of discretized time, imposing two known values $u_0$ and $u_1$ for $u$ at $s=s_0$ and $s=s_1$ and two known values $\dot{u}_0$ and $\dot{u}_1$ for $\dot{u}$ at $s=s_0$ and $s=s_1$ respectively. The amount of total computations for this method depends on the number of streamlines, on the number of time values considered, and on the spatial step size (that is, degree of discretization for the independent variable interval $[s_0,s_1]$) of the algorithm used for the resolution of a single equation. In addition, for complex geometries or turbulent flows, there is the rough-resolution of the classical Navier-Stokes system and then the interpolations of the streamlines so found.\\
Using an hardware environment for numerical computations, we have tried to limit the total time spent on this work and the amount of total memory required for handling the data, so that one can use typical workstations usually devoted to technical and scientific computing. In particular, we have used a system with 2 dual-core processors (3.2 GHz as clock frequency) and 4 GByte of central memory. We have prefered the software {\it Mathematica}, version 5.2, for the general numerics and its internal function {\ttfamily NDSolve} for the numerical resolution of differential equations.\\

The function {\ttfamily NDSolve} (see \cite{ndsolve} for details) has the option {\ttfamily Method} for specify the desired numerical scheme and the option {\ttfamily StartingStepSize} for specify the initial grid step for the independent variable. This step is then automatically adapted by an internal procedure if necessary. We have chosen for {\ttfamily Method} the {\it Adams-Moulton} or {\it Backward Differentiation} schema for their reliability when used with nonlinearity as $u \dot{u}$ type (see e.g. \cite{shampine} for technical informations and \cite{coombes} for some examples). The application of {\ttfamily NDSolve} required the specification of the initial values $u$ and $\dot{u}$ respectively for velocity and acceleration along a streamline. The numerical resolution of associated algebraic systems is then executed using advanced methods based on Lapack technology (see \cite{lapack}).\\
Version 5.2 of the software {\it Mathematica} provides support for multicore computations, spreading a single massive operation into a suitable number of parallel threads among the available processors (see \cite{wolfram}). In particular the numerical linear algebra operations seem to have the major benefits from multicore support. Therefore we have tested these features for the matrix-vector products for the computation of all the coefficients of the cubics in the case of complex shape streamlines, and for the numerical resolution of the nonlinear differential equations, where the use of implicit discretization schema involves matrix computation and resolution of many systems of algebraic equations. Also, {\it Mathematica} has the function {\ttfamily BlockMatrix} of the package {\ttfamily MatrixManipulation} for construction and manipulation of block diagonal matrix {\bf A} from block {\bf T}, and the multiplication of sparse-dense matrices like {\bf A}{\bf b} is highly optimized (see \cite{wolframMatrices}).\\

\noindent In Mathematica all the calculations are executed by a central logical unit called {\it kernel}, which distributes parallel threads of a single computation among the processors or cores available. So we have measured the performances of a single kernel spreading computations on 4 cores (2 physical processors), simulating a shared-memory system. As many kernels can run indipendently on the same machine, we have also tested the performances of 2 kernels running computations on 2 physical processors, but in this case using suitable procedures we have previously equi-distributed, among the two units, the data to handle. The two kernels have independently addressed their own memory allocation, simulating a distributed-memory system.\\
For comparison, we have executed the same computations on a single one-core processor (3.2 GHz) machine with 4 GByte of memory.

\section{About performances}
\label{performances}

The Mathematica internal function {\ttfamily AbsoluteTiming} gives the total elapsed time spent by the kernel for the execution of a single task. Therefore this time amount is the sum of the time spent by CPUs and the time spent by program (and operating system) for handling data among threads. For example, the command {\ttfamily AbsoluteTiming[A*b;]} gives the total elapsed time of execution of all the necessary operations for the calculation of the product {\ttfamily c=A*b}, without front-end visualization of the result (note the use of {\ttfamily ;} in the expression), but with product vector {\ttfamily c} stored into memory. In this way, if a task is spreaded by the kernel among parallel threads, the output given by this built-in function is the effective total time spent for completion of calculations, and it's not the sum of the single times spent by single threads.\\

\begin{table}
\label{tableMatrix}
\begin{center}
\caption{\scriptsize{Speedup for computing cubics coefficients, $100$ couples of points for each of $M$ streamlines. Right column shows speedup results in the case of 2 kernels with 1 dual processor for each of them (distributed programs), while central column shows results in the case of 1 kernel with 2 dual processors (multicore program). Multicore single program presents better performances when amount of data increases.}}

\begin{tabular}{l l l}
$M$ & \textbf{4 cores} & \textbf{2 $\times$ (2 cores)}\\
\scriptsize{1} & \scriptsize{0.7} & \scriptsize{0.35}\\
\scriptsize{10} & \scriptsize{0.9} & \scriptsize{0.6}\\
\scriptsize{50} & \scriptsize{1.3} & \scriptsize{0.95}\\
\scriptsize{100} & \scriptsize{1.85} & \scriptsize{1.1}\\
\scriptsize{200} & \scriptsize{2.7} & \scriptsize{1.3}\\
\scriptsize{300} & \scriptsize{3.15} & \scriptsize{1.5}
\end{tabular}
\end{center}
\end{table}

\noindent For an evaluation of the performances of the multicore technology, we consider the notion of {\it speedup} $S=\frac{T_0}{T}$, that is the ratio between the elapsed time spent by a single processor session and the elapsed time spent by a multicore session for the execution of the same task.\\

\noindent Now we discuss the results obtained for interpolation of streamlines in complex geometry (fig. \ref{suction}) and for resolution of corresponding flow differential equations.\\
In Table 1 we report the performances registered for the computation of $M$ streamlines interpolated by cubics, in the case of $100$ couples of points for each of them. Note the good performance of the multicore technology on single kernel case when one increases the number of streamlines, while the other case offers a slower speedup due to time spent in spreading and collecting data among the two kernels. The efficiency shown by the single kernel method when the size of the numerical parallel task increases is an experimental proof of the so called {\it Gustafson generalization} of Amdahl Law (see e.g. \cite{pacheco}).\\

\noindent When solving differential equations, if $N$ is the number of couples of points used for interpolation, the number of ODEs for every velocity component is $N$ too, because the cubics coefficients are not the same for two consecutive pairs. So the computed numerical value of the solution at the end point of the $k$-th pair is the initial value for computing the solution of the $(k+1)$-th pair. In this way, using a built-in function like {\ttfamily NDSolve}, the numerical computation of the solution along a streamline becomes an essentialy serial task. For evaluating the effect of multicore Mathematica technology in the case of {\it vectorial} operations, we have tabulated all the coefficients and initial values in a $N$ dimensional array {\ttfamily inputs} and applied {\ttfamily NDSolve} to it using a vectorial single instruction like {\ttfamily Map[modifiedNDSolve,inputs]}, being {\ttfamily Map} a built-in optimized iterator and {\ttfamily modifiedNDSolve} a suitable modified version of {\ttfamily NDSolve} which accepts as input the cubics coefficients. This algorithm is mapped to the array of all the $M$ streamlines, assigning at each step as initial values for $u$ and $\dot{u}$ the end values computed in the previous step. In particular, for $N=1$ we have the case of a priori assigned streamlines family, as that of the flow in gears pump.\\
\begin{table}
\label{tableNDSolve}
\begin{center}
\caption{\scriptsize{Speedup for computing velocity field, $N=100$ pairs of points for each of $M$ streamlines. Right column shows speedup results in the case of 2 kernels with 1 dual processor for each of them (distributed programs), while central column shows results in the case of 1 kernel with 2 dual processors (multicore program). Multicore single program has better performances than distributed parallel programs.}}

\begin{tabular}{l l l}
$M$ & \textbf{4 cores} & \textbf{2 $\times$ (2 cores)}\\
\scriptsize{1} & \scriptsize{0.55} & \scriptsize{0.4}\\
\scriptsize{10} & \scriptsize{0.8} & \scriptsize{0.6}\\
\scriptsize{50} & \scriptsize{1.1} & \scriptsize{0.75}\\
\scriptsize{100} & \scriptsize{1.4} & \scriptsize{0.9}\\
\scriptsize{200} & \scriptsize{1.85} & \scriptsize{1.1}\\
\scriptsize{300} & \scriptsize{2.3} & \scriptsize{Out of memory}
\end{tabular}
\end{center}
\end{table}

\noindent Table 2 reports the performances for the numerical resolution of ODEs. The speedups are smaller than those of the previous case, probably for the greater percentage of not parallelizable operations required by the used algorithm, but again the performances are better for the multithreads single program. Note that the virtual distributed-memory environment fails when the size of the problem reaches a critical value, while the virtual shared-memory one can complete the task, probably for a better management of memory and centrally cached informations. Note also that the elapsed time necessary for a full calculation of velocity field by a single kernel on a single processor has a range, depending on $M$, variable from few seconds to tens of minutes, therefore for great sizes or large physical domains the speedup of a multicore computation could help in a significant way.

\section{Conclusions}
\label{conclusions}

From the discussion about the used algorithms and registered performances, we can deduce the following conclusions on the parallel resolution, in a multicore environment, of blocks of ordinary differential equations:

\begin{figure}
\label{turbulence}
\begin{center}
\includegraphics[width=7cm]{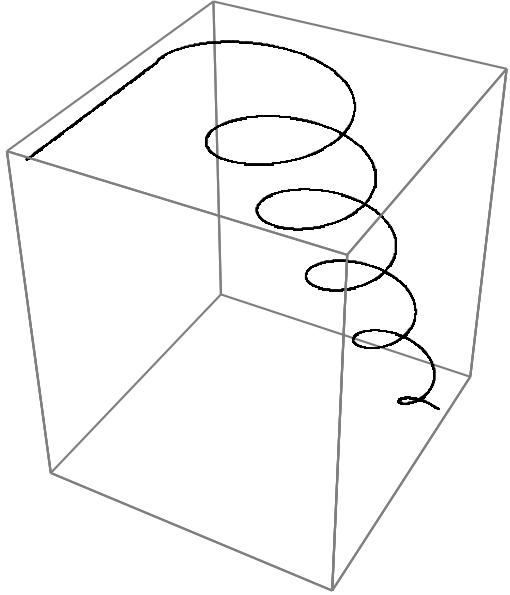}
\caption{\scriptsize{A turbulent flow streamline in the phase of oil suction. The parameterized equation is of type $s \longmapsto ([a_1+a_2s]sin(s)+a_3s,[a_1+a_2s]cos(s),a_4s)$.}}
\end{center}
\end{figure}

\begin{itemize}
\item the multicore technology gives better speedups in the case of a shared-memory environment with a single central multithreads program which executes all the computations;
\item on a multicore machine, from the last line of table (\ref{tableNDSolve}) we deduce that shared-memory single kernel resolution can execute the computations on a greater amount of data than distributed-memory multi-kernel resolution, which requires a greater allocation of memory; probably this fact depends on the use, made by shared-memory method, of a central caching of temporary data instead of a distributed one;
\item mathematical computations on large amount of scientific or technical data can take advantage from the use of numerical software specifically designed for a central multicore architecture, as already shown in previous very particular experiences (see e.g. \cite{hasslacher});
\item the streamlines algorithm and a multicore computation can help to speedup the search, in a suitable set of possible flows, of particular solutions of physical and technical interest; for example, in the case of a small gradient of pressure we have numerically tested, from a set of possible pressure fields, the existence of a solution of Navier-Stokes equations along a given turbulence streamline like that of Fig. (5), which has importance for the prediction of cavitation in pump.\\
\end{itemize} 

\noindent \textbf{Acknowledgment}\\
The author wish to thank Ing. \textit{Giuseppe Toniato}, Riello Burners's Research \& Developement Dept. Director, for his help and suggestions.\\
This work has been executed as section of a research project, supported by Italian Ministry for Scientific Research and University, on optimized flows in combustion chambers and pumps, and based on a collaboration with Torino Politecnico (Italy).



\end{document}